\input amstex
\documentstyle{amsppt}
\magnification=\magstep1
\hsize=5in
\vsize=7.3in
\TagsOnRight
\topmatter
\title Algebraic orbifold conformal field theories
\endtitle
\author Feng  Xu \endauthor
 
\address{Department of Mathematics, University of Oklahoma, 601 Elm Ave,
Room 423, Norman, OK 73019}
\endaddress
\email{xufeng\@ math.ou.edu}
\endemail
\abstract{     
We formulate the unitary rational orbifold conformal field theories
in the algebraic quantum field theory framework. 
Under general conditions, we show that the orbifold of a given
unitary rational  conformal field theories generates a unitary
modular category. Many new unitary
modular categories are obtained. We also 
show that the irreducible representations of orbifolds of rank one lattice
vertex operator algebras give rise to unitary modular categories 
and determine the corresponding modular matrices, which has been
conjectured for some time. 
}
\endabstract
\endtopmatter
\document

\heading \S1. Introduction \endheading 
Cosets and orbifolds are two methods of producing new two 
dimensional conformal field theories from given ones (cf. [MS]). 
In [X2, 3,4, 5], unitary coset conformal field theories are formulated
in the  algebraic quantum field theory framework and such
a formulation is used to solve many  questions beyond the
reach of other approaches.  The main purpose of this
paper is to formulate  unitary orbifold conformal field theories
in the same framework, and to give some
applications of this formualtion. \par
There is another approach to conformal field theories by using the
theory of vertex operator algebras (cf. [B], [FLM]). In the case 
of orbifolds this has been studied for example in [DM], [DN]. 
While there are various advantages to these different approaches, 
our main results Th. 2.6 and Th. 3.6 have not been obtained previously by
other methods. \par
Under general conditions as specified in Th. 2.6, 
we show that the orbifold of  a given
unitary rational  conformal field theories generates a unitary
modular category. Th. 2.6 gives a large family of new unitary modular
categories which can be found in \S3.4 and \S4. 
The proof of  Th. 2.6 is obtained by using the results of [KLM] and
Prop. 2.5. The proof of Prop. 2.5 is based on Galois theory for von Neumann
algebras (cf. [ILP] and references therein). \par
As an application of our general theory, 
we  show in Th. 3.6  that the irreducible 
representations of orbifolds of rank one lattice
vertex operator algebras give rise to unitary modular categories  
(hence a unitary three dimensional topological quantum field theory,
cf. [Tu])
and determine the corresponding modular 
matrices. More precisely, the simple objects of the modular categories
are in one to one correspondence with  
the irreducible 
representations of these vertex operator algebras, which were classified in
[DN]. These simple objects and the modular matrices first appeared as
examples in [DVVV] based on certain heuristic arguments, and these 
examples can be clearly interpreted as a conjecture on the existence of
certain unitary modular categories with the same modular matrices.
Th. 3.6 thus confirms this conjecture. \par
We will describe the content of this paper in more details. 
\S2.1 and \S2.2 are two preliminary sections. We include
these sections partly to set up notations.
\S2.1 is  a section on sectors following  [L3] and [L4]. 
In \S2.2 the concepts of
an irreducible conformal precosheaf $\Cal A$ and 
its covariant representations as in [GL] are introduced. 
The definition  of modular matrices in [Reh] is given. Note that 
these modular matrices are very different from the definition of 
modular matrices in Chap. 13 of [Kac] even though they coincide
in all known examples. We shall refer to the modular matrices
as defined in \S2.2 as {\it genus 0 modular matrices}. We also note that
genus 0 modular matrices determine the fusion rules by definition
(cf. (7) of \S2.2). \par
In \S2.3 we first define what we call the proper action of 
a finite group $G$  on  $\Cal A$ (Definition 2.1).
We show in Prop. 2.1 
that if a finite group $G$ acts properly on  $\Cal A$, then the 
subset of  $\Cal A$ which is invariant pointwise under the action of $G$ 
gives rise to an irreducible conformal precosheaf denoted by $\Cal A^G$.
$\Cal A^G$ is called {\it the orbifold of  $\Cal A$ with respect to $G$}. 
Many
questions concerning orbifold conformal field theories can be answered
by studying covariant representations of  $\Cal A^G$.\par
In Lemma 2.2 we collect some of the results of [X1] (also cf. [BE1])
which can be applied to our setting. In Lemma 2.3 we show that 
the proper action as defined in Definition 2.1 is always outer. \par 
In \S2.4   we first recall
the definition of absolute rationality, or $\mu$-rationality of [KLM].
Prop. 2.4 concerning the 
calculation of $\mu$-indes is essentially Prop. 21 of [KLM] 
(the ideas of the proof 
appeared in \S3 of [X6]). Prop. 2.5 shows that ${\Cal A}_G$ is always
strongly additive if ${\Cal A}$ is split and strongly additive (cf. \S2.4
for the definitions).  
We use the results of [KLM] and Prop. 2.5 to show in Th. 2.6 that if 
 $\Cal A$ is $\mu$-rational, and $G$ is a finite group acting
properly on  $\Cal A$, then $\Cal A^G$ is  $\mu$-rational (with
a formula for its $\mu$ index)  and
gives rise to a unitary modular category.
In Lemma 2.7 we collect some identities from [X4] and [BEK2] which
will be used in \S3. \par
In \S3.1 we recall from [KM] the branching rules of inclusions 
$$
Spin(M)_2\subset SU(M)_1, 
Spin(M)_2\subset Spin(M)_1 \times Spin(M)_1 \subset Spin(2M)_1,
$$
where $M=2l$ is an even positive integer, and the
numbers in the subscripts are the levels of the 
positive energy representations of the corresponding 
loop groups (cf. Chap. 9 of [PS]). 
By using these branching
rules, we show in Lemma 3.1 that ${\Cal A}_{Spin(M)_2}$ 
(cf. \S2.2 for definitions)
is the orbifold of
${\Cal A}_{SU(M)_1}$ with respect to a natural ${\Bbb Z}_2$ action (complex
conjugation), and in Lemma 3.2 that the coset  
${\Cal A}_{Spin(2M)_1/Spin(M)_2}
$ is ${\Cal A}_{U(1)_{2l}}$. \par
In \S3.2 Prop. 3.3
we calculate the genus $0$ modular matrices of  ${\Cal A}_{Spin(M)_2}$ by using
Lemma 3.1 and results of [X4] and [BEK2]. One may obtain these 
results from [TL] but we believe that the ideas in \S3.2 may find
applications elsewhere. In \S3.3 Prop. 3.4 we determine the genus $0$ 
modular matrices
of the coset $ {\Cal A}_{Spin(M)_1 \times Spin(M)_1/Spin(M)_2}$ by using
the method in \S4 of [X2] and Prop. 3.3. In \S3.4  Lemma 3.5
we show that the coset 
$ {\Cal A}_{Spin(M)_1 \times Spin(M)_1/Spin(M)_2}$ is the same 
as the orbifold of ${\Cal A}_{U(1)_{2l}}$ with respect to 
a natural ${\Bbb Z}_2$ action. By comparing with [DN], the reader may
recognize that this ${\Bbb Z}_2$ action on ${\Cal A}_{U(1)_{2l}}$ corresponds
to the $-1$ isometry on the rank one lattice vertex operator algebras.\par
By using Lemma 3.5 and Prop. 3.4, we prove  Th. 3.6. \par
The reader may wonder why we formulate the orbifold of ${\Cal A}_{U(1)_{2l}}$ 
with respect to the natural ${\Bbb Z}_2$ action as a coset instead of
considering such an orbifold directly as in [DN]. The reason is that 
even though  one can use the same ideas of \S3.2 to set up a system of
equations of genus $0$ modular matrices for the orbifold,  these equations
are not sufficient to determine genus $0$  modular matrices. 
One needs to know the conformal dimensions or univalences (cf. \S2 of [GL]) 
of ``twisted representations'' (cf. The end of \S2.4 for the definition). 
The ``twisted representations'' in [DN]
is defined algebrically. To show that 
these ``twisted representations'' in [DN] give
rise to covariant representations of our orbifold, one needs 
to study  the analytical
properties of   the twisted vertex 
operators in [DN] which is not trivial if one tries this directly.\par
In \S4 we discuss more examples of orbifolds and questions.

\heading \S2.  Orbifold CFT from Algebraic QFT point of view \endheading
\subheading{\S2.1 Sectors}
Let $M$ be a properly infinite factor
and  $\text{\rm End}(M)$ the semigroup of
 unit preserving endomorphisms of $M$.  In this paper $M$ will always
be the unique hyperfinite $III_1$ factors.
Let $\text{\rm Sect}(M)$ denote the quotient of $\text{\rm End}(M)$ modulo
unitary equivalence in $M$. We  denote by $[\rho]$ the image of
$\rho \in \text{\rm End}(M)$ in  $\text{\rm Sect}(M)$.\par
 It follows from
\cite{L3} and \cite{L4} that $\text{\rm Sect}(M)$, with $M$ a properly
infinite  von Neumann algebra, is endowed
with a natural involution $\theta \rightarrow \bar \theta $  ;
moreover,  $\text{\rm Sect}(M)$ is
 a semiring.\par
Let $\rho \in \text{\rm End}(M)$ be a normal
faithful conditional expectation
$\epsilon:
M\rightarrow \rho(M)$.  We define a number $d_\epsilon$ (possibly
$\infty$) by:
$$
d_\epsilon^{-2} :=\text{\rm Max} \{ \lambda \in [0, +\infty)|
\epsilon (m_+) \geq \lambda m_+, \forall m_+ \in M_+
\}$$ (cf. [PP]).\par
 We define
$$
d = \text{\rm Min}_\epsilon \{ d_\epsilon |  d_\epsilon < \infty \}.
$$   $d$ is called the statistical dimension of  $\rho$. It is clear
from the definition that  the statistical dimension  of  $\rho$ depends only
on the unitary equivalence classes  of  $\rho$.
The properties of the statistical dimension can be found in
[L1], [L3] and  [L4].\par
Denote by $\text{\rm Sect}_0(M)$ those elements of
$\text{\rm Sect}(M)$ with finite statistical dimensions.
For $\lambda $, $\mu \in \text{\rm Sect}_0(M)$, let
$\text{\rm Hom}(\lambda , \mu )$ denote the space of intertwiners from
$\lambda $ to $\mu $, i.e. $a\in \text{\rm Hom}(\lambda , \mu )$ iff
$a \lambda (x) = \mu (x) a $ for any $x \in M$.
$\text{\rm Hom}(\lambda , \mu )$  is a finite dimensional vector
space and we use $\langle  \lambda , \mu \rangle$ to denote
the dimension of this space.  $\langle  \lambda , \mu \rangle$
depends
only on $[\lambda ]$ and $[\mu ]$. Moreover we have
$\langle \nu \lambda , \mu \rangle =
\langle \lambda , \bar \nu \mu \rangle $,
$\langle \nu \lambda , \mu \rangle
= \langle \nu , \mu \bar \lambda \rangle $ which follows from Frobenius
duality (See \cite{L2} or \cite{Y}).  We will also use the following
notation: if $\mu $ is a subsector of $\lambda $, we will write as
$\mu \prec \lambda $  or $\lambda \succ \mu $.  A sector
is said to be irreducible if it has only one subsector. \par
\vskip .1in               
\subheading 
{\S2.2 The irreducible conformal precosheaf  and its representations}
In this section we recall the notion of irreducible conformal precosheaf
and its covariant representations as described in [GL].\par
By an {\it interval} we shall always mean an open connected subset $I$
of $S^1$ such that $I$ and the interior $I' $ of its complement are
non-empty.  We shall denote by  ${\Cal I}$ the set of intervals in $S^1$.
We shall denote by  $PSL(2, {\bold R})$ the group of
  conformal transformations on the complex plane
that preserve the orientation and leave the unit circle $S^1$ globally
invariant.  Denote by ${\bold G}$
the universal covering group of $PSL(2, {\bold R})$.  Notice that  ${\bold G}$
is a simple Lie group and has a natural action on the  unit circle
$S^1$. \par 
Denote by  $R(\vartheta )$  the (lifting to ${\bold G}$ of the) rotation by
an angle $\vartheta $. This one-parameter subgroup of
${\bold G}$ will be referred to as rotation group (denoted by Rot)
in the following.
We may associate a
one-parameter group with any interval $I$ in the following way.
Let $I_1$ be the upper
semi-circle, i.e. the interval
$\{e^{i\vartheta }, \vartheta \in (0, \pi )\}$.
 By using the Cayley transform
$C:S^1 \rightarrow {\bold R} \cup \{\infty \}$ given by
$z\rightarrow -i(z-1)(z +1)^{-1}$,
we may identify  $I_1$
with the positive real line ${\bold R}_+$. Then we consider
the one-parameter group $\Lambda _{I_1}(s)$ of diffeomorphisms of
$S^1$  such that
$$
C\Lambda _{I_1} (s) C^{-1} x = e^s x \, ,
\quad  s, x\in {\bold R} \, .  
C\Lambda _{I_1} (s) C^{-1} x = e^s x \, ,
\quad  s, x\in {\bold R} \, .
$$
We also associate with $I_1$ the reflection $r_{I_1}$ given by
$$
r_{I_1}z = \bar z
$$
where $\bar z$ is the complex conjugate of $z$.  It follows from
the definition  that
$\Lambda _{I_1}$ restricts to an orientation preserving diffeomorphisms of
$I_1$, $r_{I_1}$ restricts to an orientation reversing diffeomorphism of
$I_1$ onto $I_1^\prime $. \par
Then, if $I$ is an interval and we choose $g\in {\bold G}$ such that
$I=gI_1$ we may set
$$
\Lambda _I = g\Lambda _{I_1}g^{-1}\, ,\qquad
r_I = gr_{I_1}g^{-1}\, .
$$      
Let $r$ be an orientation reversing isometry of $S^1$ with
$r^2 = 1$ (e.g. $r_{I_1}$).  The action of $r$ on $PSL(2, {\bold R})$ by
conjugation lifts to an action $\sigma _r$ on ${\bold G}$, therefore we
may consider the semidirect product of
${\bold G}\times _{\sigma _r}{\bold Z}_2$.   Since
${\bold G}\times _{\sigma _r}{\bold Z}_2$ is a covering of the group
generated by $PSL(2, {\bold R})$ and $r$,
${\bold G}\times _{\sigma _r}{\bold Z}_2$ acts on $S^1$. We call
(anti-)unitary a representation $U$ of
${\bold G}\times _{\sigma _r}{\bold Z}_2$ by operators on ${\Cal H}$ such
that $U(g)$ is unitary, resp. antiunitary, when $g$ is orientation
preserving, resp. orientation reversing. \par
Now we are ready to define a  conformal
precosheaf. \par
An irreducible
conformal precosheaf ${\Cal A}$ of von Neumann
algebras on the intervals of $S^1$ 
is a map
$$
I\rightarrow {\Cal A}(I)
$$
from ${\Cal I}$ to the von Neumann algebras on a separable Hilbert space
${\Cal H}$ that verifies the following properties:
\vskip .1in
\noindent
{\bf A. Isotony}.  If $I_1$, $I_2$ are intervals and
$I_1 \subset I_2$, then
$$
{\Cal A}(I_1) \subset {\Cal A}(I_2)\, .
$$
 
\vskip .1in
\noindent
{\bf B. Conformal invariance}.  There is a nontrivial unitary
representation $U$ of
${\bold G}$  on         
${\Cal H}$ such that
$$
U(g){\Cal A}(I)U(g)^* = {\Cal A}(gI)\, , \qquad
g\in {\bold G}, \quad I\in {\Cal I} \, .
$$
\vskip .1in
\noindent
{\bf C. Positivity of the energy}.  The generator of the rotation subgroup
$U(R(\vartheta ) )$ is positive.
 
\vskip .1in
\noindent
{\bf D.  Locality}.  If $I_0$, $I$ are disjoint intervals then
${\Cal A}(I_0)$ and $A(I)$ commute.
 
The lattice symbol $\vee $ will denote `the von Neumann algebra generated
by'.
\vskip .1in
\noindent                                               
{\bf E. Existence of the vacuum}.  There exists a unit vector
$\Omega $ (vacuum vector) which is $U({\bold G})$-invariant and cyclic for
$\vee _{I\in {\Cal I}}{\Cal A}(I)$.
\vskip .1in
\noindent
{\bf F. Irreducibility}.  The only
$U({\bold G})$-invariant vectors are the scalar multiples of $\Omega$.
 
\vskip .1in
\noindent
The term irreducibility is due to the fact (cf. Prop. 1.2 of [GL]) that
under the assumption of {\bf F} $\vee_{I\in {\Cal I}} A(I) =B({\Cal H})$.
\par       
A covariant {\it representation} $\pi $ of
${\Cal A}$ is a family of representations $\pi _I$ of the
von Neumann algebras ${\Cal A}(I)$, $I\in {\Cal I}$, on a
separable Hilbert space ${\Cal H}_\pi $ and a unitary representation
$U_\pi $ of the covering group ${\bold G}$ of $PSL(2, {\bold R})$
such that the following properties hold:
$$
\align
I\subset \bar I \Rightarrow \pi _{\bar I} \mid _{{\Cal A}(I)}
= \pi _I \quad &\text{\rm (isotony)} \\
\text{\rm ad}U_\pi (g) \cdot \pi _I = \pi _{gI}\cdot
\text{\rm ad}U(g) &\text{\rm (covariance)}\, .
\endalign
$$
A covariant representation $\pi$ is called irreducible if
 $\vee _{I\in {\Cal I}}\pi({\Cal A}(I)) = B({\Cal H}_\pi)$. By our definition
the irreducible conformal precosheaf
 is in fact an irreducible representation  
of itself
and we will call this representation the {\it vacuum representation}. \par
Let $H$ be a simply connected simply-laced compact Lie group. By Th. 3.2
of [FG], 
the vacuum positive energy representation of the loop group
$LH$ (cf. [PS]) at level $k$ 
gives rise to an irreducible conformal precosheaf 
denoted by {\it ${\Cal A}_{H_k}$}. By Th. 3.3 of [FG], every 
irreducible positive energy representation of the loop group
$LH$ at level $k$ gives rise to  an irreducible covariant representation 
of ${\Cal A}_{H_k}$. When $H_k\subset G_l$ is a connected subgroup of a simply
connected Lie group,  Prop. 2.2 in [X2] gives an irreducible
conformal precosheaf which will be denoted by {\it ${\Cal A}_{G_l/H_k}$}
and this is refered to as {\it the coset conformal precosheaf}.   
We will see such examples in \S3.\par 
Next we will recall some of the results of [Reh]  and
introduce
notations. \par                             
Let $\{[\rho_i], i\in I \}$ be a finite set  of
equivalence classes of irreducible
covariant representations of an irreducible conformal precosheaf
with finite index . For the definitions of the conjugation and
composition of covariant representations, see \S4 of [FG] or \S2 of [GL].\par
Suppose this set is closed under
conjugation
and composition. We will denote the conjugate of $[\rho_i]$ by
$[\rho_{\bar i}]$
and identity sector by $[1]$ if no confusion arises, and let
$N_{ij}^k = \langle [\rho_i][\rho_j], [\rho_k]\rangle $. 
Here $\langle x,y\rangle$ denotes the dimension of the
space of intertwinners from $x$ to $y$ (denoted by $\text {\rm
Hom}(x,y)$) for any 
representations $x$ and $y$ (By Th. 2.3 of [GL] we don't have
to distinguish between local and global intertwinners here).  We will
denote by $\{T_e\}$ a basis of isometries in $\text {\rm
Hom}(\rho_k,\rho_i\rho_j)$.
The univalence of $\rho_i$ and the statistical dimension of
(cf. \S2  of [GL]) will be denoted by
$\omega_{\rho_i}$ and $d_{\rho_i}$ respectively. \par
Let $\phi_i$ be the unique minimal
left inverse of $\rho_i$, define:   
$$
Y_{ij}:= d_{\rho_i}  d_{\rho_j} \phi_j (\epsilon (\rho_j, \rho_i)^*
\epsilon (\rho_i, \rho_j)^*), \tag 0
$$ where $\epsilon (\rho_j, \rho_i)$ is the unitary braiding operator
 (cf. [GL] ). \par
We list two properties of $Y_{ij}$ (cf. (5.13), (5.14) of [Reh]) which
will be used in the following:
$$
\align
Y_{ij} = Y_{ji} & = Y_{i\bar j}^* = Y_{\bar i \bar j} \tag 1 \\
Y_{ij} = \sum_k N_{ij}^k \frac{\omega_i\omega_j}{\omega_k} d_{\rho_k} \tag
2
\endalign
$$
 Define
$\tilde \sigma := \sum_i d_{\rho_i}^2 \omega_{\rho_i}^{-1}$.
If the matrix $(Y_{ij})$ is invertible,
by Proposition on P.351 of [Reh] $\tilde \sigma$ satisfies
$|\tilde \sigma|^2 = \sum_i d_{\rho_i}^2$.
Suppose $\tilde \sigma= |\tilde \sigma| \exp(i x), x\in {\Bbb R}$. 
Define matrices
$$
S:= |\tilde \sigma|^{-1} Y, T:=  C Diag(\omega_{\rho_i})
\tag 3
$$ where $C:= \exp(i \frac{x}{3}).$  Then these matrices satisfy the algebra:
$$
\align
SS^{\dag} & = TT^{\dag} =id, \tag 4  \\
TSTST&= S, \tag 5 \\
S^2 =\hat{C}, T\hat{C}=\hat{C}T=T, \tag 6
\endalign
$$
where $\hat{C}_{ij} = \delta_{i\bar j}$ is the conjugation matrix. Moreover
$$
N_{ij}^k = \sum_m \frac{S_{im} S_{jm} S_{km}^*}{S_{1m}}. \tag 7
$$
(7) is known as Verlinde formula. \par
We will refer the $S,T$ matrices
as defined in  (3)  as  {\bf genus 0 modular matrices of ${\Cal A}$} since
they are constructed from the fusions rules, monodromies and minimal 
indices which can be thought as  genus 0 data associated to
a Conformal Field Theory (cf. [MS]). \par
It follows from (7) and (4) that any irreducible representation
of the commutative ring generated by $i$'s is of the form
$i\rightarrow \frac{S_{ij}}{S_{1j}}$. \par  
\subheading {\S2.3. The orbifolds}
Let ${\Cal A}$ be an irreducible conformal precosheaf on a Hilbert space
${\Cal H}$ and let $G$ be a finite group. Let $V:G\rightarrow U({\Cal H})$
be a faithful\footnotemark\footnotetext{
If $V:G\rightarrow U({\Cal H})$ is not faithful, we can take 
$G':= G/ker V$ and consider $G'$ instead.} 
unitary representation of $G$ on ${\Cal H}$.
\proclaim{ Definition 2.1} 
We say that $G$ acts properly on ${\Cal A}$ if the following conditions
are satisfied:\par
(1) For each fixed interval $I$ and each $g\in G$, 
$\alpha_g (a):=V(g)aV(g^*) \in {\Cal A}(I), \forall a\in
{\Cal A}(I)$; \par
(2) For each  $g\in G$ and $h\in {\bold G}$, $[V(g), U(h)]=0$;
moreover $V(g)\Omega = \Omega, \forall g\in G$.\par
\endproclaim
Suppose that a finite group $G$ acts properly on   ${\Cal A}$ as above.
For each interval $I$, define ${\Cal B}(I):=\{a\in {\Cal A}(I)| 
V(g)aV(g^*)=a, \forall g\in G \}$. Let $ {\Cal H}_0= \{x\in {\Cal H}|
V(g) x=x, \forall g\in G \}$ and $P_0$ the projection from ${\Cal H}$ to 
${\Cal H}_0$. Notice that $P_0$ commutes with every element of
${\Cal B}(I)$ and $U(g), \forall g\in {\bold G}$. 

Define ${\Cal A}^G(I):={\Cal B}(I)P_0$ on ${\Cal H}_0$. The unitary
representation $U$ of ${\bold G}$ on ${\Cal H}$ restricts to
an  unitary
representation (still denoted by $U$) of ${\bold G}$ on ${\Cal H}_0$.
Then:
\proclaim{2.1 Proposition} 
The map $I\in {\Cal I}\rightarrow {\Cal A}^G(I)$ on $ {\Cal H}_0$ 
together with the  unitary
representation (still denoted by $U$) of ${\bold G}$ on ${\Cal H}_0$
is an
irreducible conformal precosheaf.
\endproclaim
\demo{Proof} 
We need to check conditions ${\bold A}$ to $ {\bold F}$. ${\bold A}, 
 {\bold B},  {\bold C}, {\bold D}, {\bold F}$ are immediate consequences of 
definitions. To check ${\bold E}$, we have to show that the vacuum vector
$\Omega$ is cyclic for $\vee_{I\in {\Cal I}} {\Cal A}^G(I)$ on ${\Cal
H}_0$. Fix an interval $I$ and let $a\in {\Cal A}(I)$. Then by definition
$P_0(a\Omega)= \frac{1}{|G|} \sum_g V(g)a\Omega$. Since $ V(g)\Omega=
\Omega$, we have
$$
P_0(a\Omega)= \frac{1}{|G|} \sum_g \alpha_g(a)\Omega,
$$ and so 
$P_0(a\Omega) \in {\Cal A}^G(I)\Omega$ since 
$$
\frac{1}{|G|} \sum_g  \alpha_g(a)\in {\Cal A}^G(I).
$$ Since ${\Cal A}(I)\Omega$ is dense in ${\Cal H}$ by Reeh-Schlieder
theorem (cf. Prop. 1.1 of [GL]), 
it follows that  ${\Cal A}^G(I)\Omega$ is dense in
${\Cal H}_0$, and this shows  ${\bold E}$. 
\enddemo
\qed
\par  
The irreducible conformal precosheaf in Prop. 2.1 will be denoted by
${\Cal A}^G$ and will be called the {\it orbifold of ${\Cal A}$}
with respect to $G$. \par
The net ${\Cal B}(I)\subset {\Cal A}(I)$ 
is a standard net of inclusions (cf. [LR])
with condtional expectation $\epsilon$ defined by 
$$
\epsilon:=\frac{1}{|G|} \sum_g \alpha_g(a), \forall a\in {\Cal A}(I).
$$ 
Note that $\epsilon$ has finite index. We can therefore apply the theory
in [X1] to this setting (also cf. [BE1-2]). 
Fix an interval $I$ and let $M:= {\Cal A}^G(I)$. 
Denote by $i$ 
(resp. $\lambda$) the irreducible covariant representations of ${\Cal A}^G$
(resp. ${\Cal A}$) with finite index, and $\gamma$ the restriction of 
the vacuum representation of $ {\Cal A}$ to ${\Cal A}^G$. Denote by
$b_{i\alpha}\in {\Bbb Z}$ 
the multiplicty of representation $\lambda$ which appears
in the restriction of  representation $i$ when restricting from 
$ {\Cal A}$ to ${\Cal A}^G$. $b_{i\alpha}$ is also known as the branching
rules.   
Then there are  maps 
$$
\lambda\rightarrow a_{\lambda}, i\rightarrow \sigma_i
$$
where $a_{\lambda}, \sigma_i\in \text{\rm End} (M)$
with some remarkable properties. We list some of the properties in
the following:  
\proclaim{Lemma 2.2}
The map $
\lambda\rightarrow a_{\lambda}
$ as a map of sectors is a ring homomorphism. The map 
$i\rightarrow \sigma_i$ as a map of sectors is a ring isomorphism. 
Moreover
$$ \langle a_{\lambda},  a_{\mu} \rangle =
\langle \lambda, \mu \gamma \rangle, b_{i\lambda}= \langle \sigma_i,
  a_{\lambda} \rangle.
$$
\endproclaim
For the proof of Lemma 2.2, we refer the reader to \S3 of [X1] or [BE1]. 
Note that the proof in \S3 of [X1] is given in the case of conformal
inclusions, but the same proof applies verbatim.
We also note that in the proof of above equations 
(cf. [X1] or [BE1]) strong additivity on ${\Cal B}(I)$
is assumed to ensure the equivalence of local and global intertwinners
. More precisely the equivalence of local and global intertwinners
is used in proving Braiding-Fusion equations. The strong additivity
assumption is not necessary under the condition of conformal invariance
and finite index since under these conditions 
the equivalence of local and global intertwinners has been proved in
Th. 2.3 of [GL].  \par
As noticed in [BE1], if $\lambda=\mu$ is the identity sector, then from
above we have:
$\langle id,  \gamma \rangle =1$, which implies that the 
inclusion ${\Cal B}(I)\subset {\Cal A}(I)$ is irreducible, i.e.,
$ {\Cal B(I)}'\cap {\Cal A}(I) \simeq {\Bbb C}$. 
We record this result in the following:
\proclaim{Lemma 2.3}
If $G$ acts properly on ${\Cal A}$ as defined above, then for each
$I$, the action of $G$ on ${\Cal A}(I)$ is outer, i.e.,
$ {\Cal B(I)}'\cap {\Cal A}(I) \simeq {\Bbb C} $ where 
${\Cal B}(I)$ is the fixed point subalgebra
of ${\Cal A}(I)$ under the action of $G$.
\endproclaim
Since the action of $G$ on ${\Cal A}(I)$ is outer, there exists a unique
normal faithful conditional expectation $\epsilon: {\Cal A}(I)\rightarrow
{\Cal B}(I)$ given by
$$
\epsilon(x)=  \frac{1}{|G|}\sum_{k\in K} \alpha_g(x).
$$
The index of $\epsilon$ is $|G|$ and this is the index of the 
inclusion ${\Cal B}(I)\subset  {\Cal A}(I)$. There exists an isometry
$v$ in $ {\Cal A}(I)$ (cf. \S2 of [LR] or references therein) such that
$$
 \epsilon(xv^*)v= v^*\tilde \epsilon(vx)= \frac{1}{|G|} x,
\epsilon(vv^*) =  \frac{1}{|G|}, \alpha_g(V)^* \alpha_h(V)= \delta_{g,h}
id.
$$
Moreover the restriction of the vacuum 
representation of $ {\Cal A}$ to ${\Cal A}^G$, denoted by $\gamma$ in Lemma 
2.2, decomposes as $\gamma= \sum_{\lambda} b_{1\lambda}\lambda$, and
$b_{1\lambda}= d_\lambda$ where $d_\lambda$ is the statistical dimension
of $\lambda$. We have 
$$
\sum_\lambda b_{1\lambda}d_\lambda= \sum_\lambda (d_\lambda)^2 = |G|.
$$
\subheading{\S2.4. Absolute rationality }
We first recall some definitions from [KLM] .
As in [GL] by an interval of the circle we mean an open connected
proper subset of the circle. If $I$ is such an interval then
$I'$ will denote the interior of the complement of $I$ in the circle.
We will denote by ${\Cal I}$ the set of such intervals.
Let $I_1, I_2\in {\Cal I}$. We say that $I_1, I_2$ are disjoint if
$\bar I_1\cap \bar I_2=\emptyset$, where $\bar I$
is the closure of $I$ in $S^1$.  
When $I_1, I_2$ are disjoint, $I_1\cup I_2$
is called a 1-disconnected interval in [X6].  
Denote by ${\Cal I}_2$ the set of unions of disjoint 2 elements
in ${\Cal I}$. Let ${\Cal A}$ be an irreducible conformal precosheaf
as in \S2.1. For $E=I_1\cup I_2\in{\Cal I}_2$, let
$I_3\cup I_4$ be the interior of the complement of $I_1\cup I_2$ in 
$S^1$ where $I_3, I_4$ are disjoint intervals. 
Let 
$$
{\Cal A}(E):= A(I_1)\vee A(I_2), 
\hat {\Cal A}(E):= (A(I_3)\vee A(I_4))'.
$$ Note that ${\Cal A}(E) \subset \hat {\Cal A}(E)$.
Recall that a net ${\Cal A}$ is {\it split} if ${\Cal A}(I_1)\vee
{\Cal A}(I_2)$ is naturally isomorphic to the tensor product of
von Neumann algebras ${\Cal A}(I_1)\otimes
{\Cal A}(I_2)$ for any disjoint intervals $I_1, I_2\in {\Cal I}$.
${\Cal A}$ is {\it strongly additive} if ${\Cal A}(I_1)\vee
{\Cal A}(I_2)= {\Cal A}(I)$ where $I_1\cup I_2$ is obtained
by removing an interior point from $I$.
\proclaim{Definition 2.2 (Absolute rationality of [KLM])}
${\Cal A}$ is said to be absolute rational, or $\mu$-rational, if
${\Cal A}$ is split, strongly additive, and 
the index $[\hat {\Cal A}(E): {\Cal A}(E)]$ is finite for some
$E\in {\Cal I}_2$ . The value of the index
$[\hat {\Cal A}(E): {\Cal A}(E)]$ (it is independent of 
$E$ by Prop. 5 of [KLM]) is denoted by $\mu_{{\Cal A}}$
and is called the $\mu$-index of ${\Cal A}$. If 
the index $[\hat {\Cal A}(E): {\Cal A}(E)]$ is infinity for some
$E\in {\Cal I}_2$, we define the $\mu$-index of ${\Cal A}$ to be
infinity.
\endproclaim
\par
We have:
\proclaim{Proposition 2.4}
Suppose ${\Cal B}\subset {\Cal A}$ is a standard net of inclusions
as defined in 3.1 of [LR]. Let $E\in {\Cal I}_2$.
If  ${\Cal A}\subset {\Cal B}$ has finite index denoted by 
$[{\Cal A}:{\Cal B}]$ and ${\Cal A}$ and ${\Cal B}$ are split, 
then
$$
[\hat {\Cal B}(E):{\Cal B}(E)] = [{\Cal A}:{\Cal B}]^2 
[\hat {\Cal A}(E):{\Cal A}(E)].
$$ 
\endproclaim
Prop. 2.4 is essentially Prop.21 of [KLM] except that we do not assume that
${\Cal A}$ and ${\Cal B}$ are $\mu$-rational. But the proof of Prop.21 of
[KLM] applies word by word (also cf. the proof of Th. 3.5 in [X6]). \par

\proclaim{Proposition 2.5}
Let ${\Cal A}$ be an irreducible conformal precosheaf and let $G$
be a finite group acting properly on ${\Cal A}$. Suppose that 
${\Cal A}$ is split and strongly additive.  Then ${\Cal A}^G$ is
also split and strongly additive.
\endproclaim
\demo{Proof:}
It follows from the definitions that ${\Cal A}^G$ is split. \par
Let $I$ be an interval, and $I_1,I_2$ are the connected components
of a set obtained from $I$ by removing an interior point of $I$. 
To show ${\Cal A}^G$ is strongly additive, it is sufficient to show that
${\Cal B}(I_1)\vee {\Cal B}(I_2) = {\Cal B}(I).$\par
Let us show that ${\Cal A}(I_1)\vee {\Cal B}(I_2)= {\Cal A}(I)$. 
First note that
$[{\Cal A}(I): {\Cal A}(I_1)\vee {\Cal B}(I_2)]<\infty$.
In fact let $I_2^{(n)}\subset I_2$ be an increasing sequence of
intervals such that $I_2^{(n)}$ have one boundary point in common with
$I_2$, $\bar I_1\cap  \bar I_2^{(n)}=\emptyset$ and
$\cup_n  I_2^{(n)}=I_2$. By the additivity of the conformal net ${\Cal A}$
(cf. \S3 of [FJ]), we have that ${\Cal A}(I_1)\vee {\Cal B}(I_2^{(n)})$
(resp. ${\Cal A}(I_1)\vee {\Cal A}(I_2^{(n)})$) are increasing
sequences of von Neumann algebras such that
$$
\vee_n {\Cal A}(I_1)\vee {\Cal B}(I_2^{(n)})= {\Cal A}(I_1)\vee {\Cal
B}(I_2),
\vee_n {\Cal A}(I_1)\vee {\Cal A}(I_2^{(n)})= {\Cal A}(I_1)\vee {\Cal
A}(I_2)= {\Cal A}(I)
$$
where we have used the assumption that ${\Cal A}$ is strongly additive.
By the splitting property we have
$$
[{\Cal A}(I_1)\vee {\Cal A}(I_2^{(n)}): {\Cal A}(I_1)\vee {\Cal
B}(I_2^{(n)})]
= [{\Cal A}(I_1)\otimes {\Cal A}(I_2^{(n)}): {\Cal A}(I_1)\otimes {\Cal
B}(I_2^{(n)})] 
= |G|. 
$$
It follows (cf. Prop. 3 of [KLM]) that
$[{\Cal A}(I): {\Cal A}(I_1)\vee {\Cal B}(I_2)]\leq |G|.$\par
So there exists a faithful normal conditional expectation
$\tilde \epsilon: {\Cal A}(I)\rightarrow {\Cal A}(I_1)\vee {\Cal B}(I_2)$.
Note that 
$$
{\Cal B}(I_2)\subset 
\tilde \epsilon({\Cal A}(I_2))\subset  {\Cal A}(I_1') \cap 
{\Cal A}(I) =  {\Cal A}(I_2)
$$ and so
$\tilde \epsilon({\Cal A}(I_2))$ in an intermediate von Neumann algebra
between ${\Cal B}(I_2)$ and ${\Cal A}(I_2)$. So (cf. [ILP] or references
therein) there exists a 
subgroup $K$ of $G$ such that $\tilde \epsilon({\Cal A}(I_2))$ is the 
pointwise fixed subalgebra of  ${\Cal A}(I_2)$ under the action of $K$.
Since ${\Cal B}(I_2)\subset {\Cal A}(I_2)$ is irreducible by Lemma 2.3, 
$\tilde \epsilon({\Cal A}(I_2))\subset {\Cal A}(I_2)$ is also irreducible
and it follows that there exists a unique  conditional expectation from
${\Cal A}(I_2)$ to $\tilde \epsilon({\Cal A}(I_2))$ and it is given by
$$
\tilde \epsilon(x_2)= \frac{1}{|K|}\sum_{k\in K} \alpha_k(x_2).
$$
Let us show that $K$ is the trivial subgroup, i.e., if $k\in K$, then
$k$ is the identity element of $G$.\par
Let $v\in {\Cal A}(I_2)$ be the isometry as in the end of \S2.2 with
$G$ replaced by $K$ such that 
$$
\tilde \epsilon(x_2v^*)v= v^*\tilde \epsilon(vx_2)= \frac{1}{|K|} x_2,
\epsilon(vv^*) =  \frac{1}{|K|}.
$$
Define a map $\tilde\gamma:{\Cal A}(I)\rightarrow {\Cal A}(I_1)\vee {\Cal
B}(I_2)$
by:
$$
\tilde\gamma(x):= |K| \tilde \epsilon(vxv^*), \forall x\in {\Cal A}(I).
$$ 
One checks easily that 
$$
\tilde\gamma(x_1)=x_1, \tilde\gamma(x_2x_2')= 
\tilde\gamma(x_2)\tilde\gamma(x_2'),
\forall x_1\in {\Cal A}(I_1), x_2,x_2'\in {\Cal A}(I_2).
$$ It follows that $\tilde\gamma(xy)=\tilde\gamma(x)\tilde\gamma(xy)$ for any
$x,y \in {\Cal A}(I)$ since $A(I)$ is generated by two
commuting subalgebras $A(I_1)$ and $A(I_2)$. \par
For any $k\in K$, define $v_k = \alpha_k(v)$ and
$$
\tilde \alpha_k(x)= v_k^* \tilde\gamma(x) v_k, \forall x\in {\Cal A}(I).  
$$
Then one checks that
$$
\align
\tilde\alpha_k(x_1)=x_1, \tilde\alpha_k(x_2x_2')= 
& \alpha_k(x_2)\alpha_k(x_2'), \\
\tilde\alpha_k(\tilde\alpha_{k^{-1}}x_2)=x_2, 
& \forall x_1\in {\Cal A}(I_1), x_2,x_2'\in {\Cal A}(I_2).
\endalign
$$ It follows that 
$$\tilde\alpha_k(xy)=\tilde\alpha_k(x)\tilde\alpha_k(xy)$$ for any
$x,y \in {\Cal A}(I)$ since $A(I)$ is generated by two
commuting subalgebras $A(I_1)$ and $A(I_2)$. One can also check similarly
$\tilde\alpha_k(\tilde\alpha_{k^{-1}}(x))=x, \forall x\in {\Cal A}(I)$. 
So $\tilde\alpha_k$ is an automorphism of ${\Cal A}(I)$. Since ${\Cal A}(I)$ is
a type $III$ factor, there exists a unitary opertor $U_k\in B(H_0)$ such
that $\tilde\alpha_k(x)= U_k x U_k^*, \forall x\in {\Cal A}(I)$. 
Since $\tilde\alpha_k(x_1)=x_1,  \forall x_1\in {\Cal A}(I_1)$, we have
$ U_k\in {\Cal A}(I_1)' ={\Cal A}(I_1')$ by Haag duality (cf. \S2 of [GL]). 
\par
If for all unitary $U'\in {\Cal A}(I')$, 
$(U_k U'\Omega, \Omega)= 0$, then $ ( U'\Omega, U_k\Omega)= 0$, and
it follows that $({\Cal A}(I') \Omega, U_k^*\Omega)= 0$. Since 
${\Cal A}(I') \Omega$ is dense in ${\Cal H}$ 
by Reeh-Schilider theorem (cf. [GL]),
it follows that  $ U_k^*\Omega=0$, which implies $ U_k^*=0$ by using
Reeh-Schilider theorem again since $U_k^*\in {\Cal A}(I_1')$. Hence there
exists a unitary $U'\in {\Cal A}(I')$ such that
$(U_k U'\Omega, \Omega)\neq 0$. Note that ${\Cal A}(I')\subset {\Cal
A}(I_1')$.
Replacing $U_k$ by $U_k U'$ if necessary, we may assume that
$(U_k \Omega, \Omega)\neq 0$. \par
Let $g_n\in {\bold G}$ be a sequence of elements such that
$g_n I_1' = I_1'$ and $g_n I_2$ is an increasing sequence intervals 
containing $I_2$, i.e.,
$I_2\subset g_n I_2\subset g_{n+1} I_2$, and
$\cup_n  g_n I_2=I_1'$ (One may take $g_n$ to be a sequence of
dilations). Applying $Ad(U(g_n))$ to the equation
$$
U_k x_2= \alpha_k(x_2) U_k
$$ and using $\alpha_k(Ad(U(g_n)) x_2)= Ad(U(g_n)) (\alpha_k (x_2))$, we get
$$
Ad(U(g_n))(U_k) Ad(U(g_n))(x_2) = \alpha_k(Ad(U(g_n)) x_2) Ad(U(g_n))(U_k).
$$
It follow that 
$$
Ad(U(g_n))(U_k) x_2^{(n)} = \alpha_k( x_2^{(n)}) Ad(U(g_n))(U_k),
\forall x_2^{(n)}\in {\Cal A}( g_n I_2).
$$
Let $U$ be a weak limit of $Ad(U(g_n))(U_k)$.  Note that
$$U\in{\Cal A}(I_1')
$$ since $Ad(U(g_n))(U_k)\in {\Cal A}(I_1')$ by 
our choice of $g_n$ .  Since
$(Ad(U(g_n))(U_k) \Omega, \Omega)= (U_k \Omega, \Omega)\neq 0$ where
we use the fact that $\Omega$ is invariant under the action of $U(g_n)$,
it follows that
$$
(U  \Omega, \Omega)= (U_k \Omega, \Omega)\neq 0,
$$ so $U\neq 0$, and we have
$$
U x_2^{(n)} = \alpha_k( x_2^{(n)}) U,
\forall x_2^{(n)}\in {\Cal A}( g_n I_2).
$$
Since $\cup_n  g_n I_2=I_1'$, $\vee_n {\Cal A}( g_n I_2)= A(I_1')$ by
the additivity of the conformal net ${\Cal A}$ (cf. \S3 of [FJ]), 
it follows that
$$
U\neq 0, U\in{\Cal A}(I_1'),  U x = \alpha_k( x) U,
\forall x\in {\Cal A}(I_1').
$$
Recall that $\alpha_k$ is an automorphism of ${\Cal A}(I_1')$ and
${\Cal A}(I_1')$ is a factor, it follows that
$$
UU^* = c id = U^*U, c\neq 0. 
$$
Change $U$ into $\frac{1}{\sqrt c} U$ if necessary, we may assume that
$U$ is unitary, and so we have
$$
\alpha_k(x)= Ad U (x), \forall x\in {\Cal A}(I_1').
$$ So $U\in B(I_1')'\cap {\Cal A}(I_1')$, and by Lemma 2.3, 
$Ad U (x)= x =\alpha_k(x) , \forall x\in {\Cal A}(I_1')$.  
It follow that $V_k x\Omega = x\Omega, \forall   x\in {\Cal A}(I_1')$,
and by Reeh-Schilder theorem $V_k=id$, so $k$ is the identity element 
in $G$. Since $k\in K$ is arbitary, we have shown that $K$ is the
trivial group. \par
So ${\Cal A}(I_1)\vee {\Cal B}(I_2)= {\Cal A}(I)$, and
$$
{\Cal B}(I)=
\epsilon ({\Cal A}(I_1)) = \epsilon ({\Cal A}(I_1)\vee {\Cal B}(I_2))
= {\Cal B}(I_1)\vee {\Cal B}(I_2).
$$
\enddemo
\qed \par
\proclaim{Theorem 2.6}
Let ${\Cal A}$ be an irreducible conformal precosheaf and let $G$
be a finite group acting properly on ${\Cal A}$. Suppose that 
${\Cal A}$ is absolutely rational or $\mu$-rational as in definition
2.2. Then:\par
(1): ${\Cal A}^G$ is absolutely rational or $\mu$-rational and  
$\mu_{{\Cal A}^G}= |G|^2 \mu_{{\Cal A}}$; \par
(2): There are only a finite number of irreducible covariant 
representations of ${\Cal A}^G$ and they give rise to a unitary modular
category as defined in II.5 of [Tu] by the consturction as given in
\S1.7 of [X5].
\endproclaim
\demo{Proof:}
Ad (1):
Note that ${\Cal B}(I):=\{a\in {\Cal A}(I)| 
V(g)aV(g^*)=a, \forall g\in G \} \subset {\Cal A}(I)$ is a standard net
with
index $|G|$ by our definitions. By Prop. 2.2, for any $E\in {\Cal I}_2$,  
$$
[\hat {\Cal A}^G(E):{\Cal A}^G(E)] = |G|^2 
[\hat {\Cal A}(E):{\Cal A}(E)] =  |G|^2 \mu_{{\Cal A}}.
$$  This proves (1). \par
Ad (2): 
This follows from Prop. 2.5, Cor. 32 of [KLM] 
and (2) of Cor. 1.7.3 of [X5].
\enddemo
\qed
\par
Suppose that ${\Cal A}$ and $G$ satisfy the assumptions of Th. 2.6. 
By Th. 2.6 and Th. 30 of [KLM],  ${\Cal A}^G$ has only finite 
finite number of irreducible representations $\lambda$ and
$$
\sum_{\lambda}d_\lambda^2 = \mu_{{\Cal A}^G}.
$$
The set of such $\lambda$'s is closed under conjugation and compositions,
and by Cor. 32 of [KLM], the $Y$-matrix as defined in \S2.1 for
${\Cal A}^G$ is non-degenerate, and we will denote the corresponding
genus $0$ modular matrices by $\dot S, \dot T$. 
Let $\dot C$ (resp. $C$) be the constant as defined in (3) of \S2.1
for ${\Cal A}^G$ (resp. ${\Cal A}$).
Denote 
the genus $0$ modular matrices of ${\Cal A}$ by $ S,  T$. 
Let $W$ be the vector space over ${\Bbb C}$ with basis given 
the irreducible subsectors of $\sigma_ia_\lambda, \forall i, \lambda$.
\proclaim{Lemma 2.7}
Under the conditions above we have:\par
(1): $\sum_\lambda b_{i\lambda} \dot S_{\lambda \nu} = 
\sum_k S_{i k} b_{k\nu}$; \par
(2): $\dot C^3= C^3$; \par
(3): The dimension of the vector space $W$ is 
$\sum_{i,\lambda}  b_{i\lambda}^2$.
\endproclaim 
\demo{Proof}
(1) and (2) follows from the proof of Prop. 3.1 in [X4] (also cf.
the remark after Prop. 3.1 in [X4]). (3) follows from Cor. 4.12 of
[BEK2].
\enddemo
\qed
\par  
An irreducible covariant representation  $\lambda$  of ${\Cal A}^G$
is called
an {\it untwisted} representation if $b_{i\lambda}\neq 0$ for some $i$.
These are
representations of ${\Cal A}^G$ which appear as subrepresentations in the
the restriction of some representation of ${\Cal A}$ to 
${\Cal A}^G$. A  representation is called {\it twisted} if it is not
untwisted.    
Note that 
$\sum_{\lambda} d_\lambda b_{i\lambda} = d_i |G|$, and
$ b_{1\lambda} =d_\lambda$. So we have
$$
\sum_{\lambda \ \text{\rm untwisted} } d_\lambda^2 \leq
\sum_i({\sum_{\lambda} d_\lambda b_{i\lambda}})^2
= |G|+ \sum_{i\neq 1} d_i^2 |G|^2
< |G|^2+ \sum_{i\neq 1} d_i^2 |G|^2 = \mu_{{\Cal A}^G}
$$ 
if $G$ is not a trivial group, where in the last $=$ we have used 
Th. 2.6. 
It follows that the set of twisted representations of 
${\Cal A}^G$ is not empty. This is very different from the case of
cosets, cf. [X4] Cor. 3.2 where it was shown that under certain
conditions there are no twisted representations for the coset.\par

\heading \S3. A class of orbifolds  \endheading 
\subheading{\S3.1 Some inclusions}
We recall from \S4 of [KM] the  branching rules of inclusions
$$
Spin(M)_2\subset SU(M)_1,
Spin(M)_2\subset Spin(M)_1 \times Spin(M)_1 \subset Spin(2M)_1,
$$
where $M=2l$ is an even positive integer and $l\geq 3$. \par
We will use $H_x$ to denote the Hilbert space of the 
positive energy representation of a loop group (cf. \S9 of [PS])
with integrable weight $x$. \par 
For the inclusion $Spin(M)_2\subset SU(M)_1$, we have
$$
\align
H_{\tilde \Lambda_0}= H_{2 \Lambda_0} + H_{2\Lambda_1},
& H_{\tilde \Lambda_l}= H_{2 \Lambda_l} + H_{2\lambda_{l-1}}; \\
H_{\tilde \Lambda_1}= H_{\Lambda_0+\Lambda_1}, 
& H_{\tilde \Lambda_{2l-1}}= H_{\Lambda_0+\Lambda_1}; \\
H_{\tilde \Lambda_{l-1}}= H_{\Lambda_{l-1}+\Lambda_1 }, 
& H_{\tilde \Lambda_{l+1}}= H_{\Lambda_{l-1}+\Lambda_1 }; \\
H_{\tilde \Lambda_i}= H_{\Lambda_i} 
&  H_{\tilde \Lambda_{2l-i}}= H_{\Lambda_i}, 2\leq i\leq l-2. \tag 1
\endalign
$$
Where $\tilde\Lambda_i, 0\leq i\leq 2l-1$ are the $2l$ integrable weights
of $LSU(M)$ at level $1$, and 
$$
2 \Lambda_0, 2\Lambda_1, 2 \Lambda_l, 2\lambda_{l-1},
\Lambda_0+\Lambda_1,\Lambda_{l-1}+\Lambda_1, \Lambda_i, 2\leq i\leq l-2
$$ are integrable weights of $LSpin(M)$ at level $2$.\par
$$
H_{\dot \Lambda}\otimes H_{\ddot \Lambda} =\sum_{\Lambda, \dot \Lambda+
\ddot \Lambda - \Lambda\in Q} H_{(\dot \Lambda,\ddot \Lambda;  \Lambda)}
\otimes H_{\Lambda} \tag 2
$$
Where $Q$ is the root lattice of $Spin (M)$ (cf. \S1.3 of [KW]), and 
$(\dot \Lambda, \ddot \Lambda)$ are the integrable weights (both at level
$1$)
of $LSpin(M)\times LSpin(M)$.\par 

The irreducible conformal precosheaf ${\Cal A}_{U(1)_{2l}}$ associated with
$LU(1)$ at level $2l$ is studied in \S3.5 of [X5]. We have
$\mu_{{\Cal A}_{U(1)_{2l}}}= 2l$, and there are exactly $2l$ irreducible
representations of ${\Cal A}_{U(1)_{2l}}$ which is labeled by integer 
$k, 0\leq k\leq
2l-1$. 
By identifying ${\Bbb R}^{2M}=(x,y)\rightarrow x+iy\in {\Bbb C}^{M}$
where $x, y$ are column vectors with $M$ real entries, we have the 
following natural inclusion  
$SU(M)_1 \times U(1)_M \subset Spin (2M)_1$ where $U(1)$ acts on
${\Bbb C}^{M}$ as a complex scalar. We have:
$$
H_{\hat \Lambda_0}= \sum_{0\leq i\leq l} H_{\tilde \Lambda_{2i}}\otimes
H_{2i} \tag 3
$$
where $H_{\hat \Lambda_0}$ is the vacuum representation of $LSpin(2M)$ at
level $1$. \par
Note that we have natural inclusion
$Spin(M)_2\subset SU(M)_1 \subset Spin (2M)_1.
$
Define $J:= (Id_M, -Id_M) \in SO(2M)$ and lift it to $Spin (2M)$. 
Note that for $A\in SU(M), JAJ= \bar A$, 
and $JAJ=A$ if $A\in Spin(M)$. 
The operator $J$ preserves $H_{\tilde \Lambda_0}$, and one checks
$$
J.x= x, \forall x\in H_{\Lambda_0}, 
J.x= -x, \forall x\in H_{2\Lambda_1}.  
$$
Moreover on $H_{\hat \Lambda_0}$ $J$ commutes with the action
of ${\bold G}$.  
It follows that $Ad J$ generates a proper ${\Bbb Z}_2$ action on
${\Cal A}_{SU(M)_1}$, and by the first equation in (1) we have:
\proclaim{Lemma 3.1}
$$
{\Cal A}_{Spin(M)_2} = {\Cal A}_{SU(M)_1}^{{\Bbb Z}_2}.
$$
\endproclaim
The next Lemma will be used in \S3.4.
\proclaim{Lemma 3.2}
${\Cal A}_{Spin(2M)_2/Spin(M)_1} = {\Cal A}_{U(1)_M}$.
\endproclaim
\demo{Proof}
From definition we have 
$${\Cal A}_{U(1)_{M}}(I) \subset {\Cal A}_{Spin(2M)_2/Spin(M)_1}(I).
$$
Since the modular automorphism group
of ${\Cal A}_{Spin(2M)_2/Spin(M)_1}(I)$ with respect to the vacuum
vector $\Omega$ is geometric (cf. \S2 of [GL]) and fixes globally
${\Cal A}_{U(1)_{M}}(I)$, by Takesaki's theorem (cf. [T])
we just have to show that 
$$\overline{{\Cal A}_{U(1)_{M}}(I)\Omega} = 
\overline{{\Cal
A}_{Spin(2M)/Spin(M)}(I)
\Omega}.$$ 
By Reeh-Schilider's theorem (cf. \S2 of [GL]) it is sufficient to show that
$$\overline{LU(1)_{M}\Omega} = 
\overline{\vee_{I\in {\Cal I}}{\Cal
A}_{Spin(2M)/Spin(M)}(I).
\Omega}.
$$
Note that by the branching rules (1) and (3) 
$$ \overline{\vee_{I\in {\Cal I}}{\Cal
A}_{Spin(2M)/Spin(M)}(I)
\Omega}\subset \overline{LU(1)_{M}\Omega}
$$ and the proof is complete. 
\enddemo
\qed \par
\subheading {\S3.2. Genus $0$ modular matrices of $Spin(M)_2$}
In this section we determine the genus $0$ 
modular matrices of  ${\Cal A}_{Spin(M)_2}$
by using Lemma 3.1 and Th. 2.6. By Th. 2.6,  ${\Cal A}_{Spin(M)_2}$ is absolute
rational with $\mu$-index equal to $4M$. By Cor.9  of [KLM], each 
irreducible representation of ${\Cal A}_{Spin(M)_2}$ has finite index. 
The following is a list of irreducible representations of 
${\Cal A}_{Spin(M)_2}$ (as before we use the integrable weights to denote the
corresponding representations):
$$
\align
& \widehat{1}:=2\Lambda_0, \widehat{j}:=2\Lambda_1, 
\widehat{\phi_l^1}:=2\Lambda_{l-1}, \widehat{\phi_l^2}=2\Lambda_{l}; \\
& \widehat{\phi_1}:=\Lambda_0+ \Lambda_1, 
\widehat{\phi_{l-1}}:=\Lambda_{1-1}+ \Lambda_{1},
\widehat{\phi_2}:=\Lambda_2,..., \widehat{\phi_{l-2}}:=\Lambda_{l-2}; \\
& \widehat{\sigma_1}:=\Lambda_0+\Lambda_{l-1}, 
\widehat{\sigma_2}:+\Lambda_0+\Lambda_{1}, 
 \widehat{\tau_2}:=\Lambda_1+\Lambda_{l-1}, 
 \widehat{\tau_1}:=\Lambda_1+\Lambda_{l} \tag 4
\endalign
$$
where we have chosen our notations so that one may easily compare with
the notations of [DVVV]. 
We will show that the above is in fact a complete list of
irreducible  representations of ${\Cal A}_{Spin(M)_2}$.\par
The center $Z$ of $Spin(M)$ is ${\Bbb Z}_4$ when $l$ is odd, and
${\Bbb Z}_2\times {\Bbb Z}_2$ when $l$ is even, and 
$Z$  acts transitively (cf. Chap. 9 of [PS]) on the set 
$2\Lambda_0, 2\Lambda_1, 2\Lambda_{l-1}, 2\Lambda_{l}$. The action
is also known as diagram automorphisms as defined in \S3.3.
Choose a path (localized on an interval $I$)
$P: [0,2\pi]\rightarrow H$ where $H\subset Spin(M)$ is the Cartan subgroup
of $Spin(M),$ such that $P(0)= id$ of $Spin(M)$, and $P(2\pi)\in Z.$
$AdP$ gives an outer action of $LSpin(M)$ if 
$P(2\pi)\neq id$.  One has (cf. Chap. 9 of [PS]):
$$ \pi_{\Lambda} (AdP.LSpin(M)) \simeq \pi_{P(0).\Lambda}(LSpin(M)),
$$ where
$P(0).\Lambda$ denotes the image of the action of $P(0)$ on 
$\Lambda$. From this one can determine the
following fusion rules (See \S3.3 for the definition of 
diagram automorphisms): 
$$
\align
\widehat{\phi_l^1}^2 =\widehat{j}, \widehat{\phi_l^1}^3=
\widehat{\phi_l^2}, \widehat{\phi_l^1}^4= \widehat{1},
& \widehat{\phi_l^1}\widehat{\sigma_1} = \widehat{\tau_2},
\text {\rm
if} \ l\in 2{\Bbb Z}+1; \\
\widehat{\phi_l^1}^2=\widehat{\phi_l^2}^2= \widehat{j}^2 =
\widehat{1},
& \widehat{\phi_l^1}\widehat{\sigma_1} = \widehat{\sigma_1},
\widehat{\phi_l^2}\widehat{\sigma_2} = \widehat{\sigma_2},
\widehat{j}\widehat{\sigma_k}= \widehat{\tau_k}, k=1,2,
\text {\rm
if} \ l\in 2{\Bbb Z}. \tag 5
\endalign
$$
We will denote by $\dot S, \dot T$ (resp. $S, T$) the genus
$0$ modular matrices of ${\Cal A}_{Spin(M)_2}$ (resp. 
${\Cal A}_{SU(M)_1}$). Note that by [W2]
$S_{\Lambda_k, \Lambda_{k'}} = \frac{1}{\sqrt{2l}}
\exp(\frac{2\pi i kk'}{2l}), 0\leq k,k'\leq 2l-1.$ \par
Consider the space $W$ as defined before Lemma 2.7.
We have $dim W= \sum_{i,\alpha} b_{i\alpha}^2 =M+2$ by (3)
of Lemma 2.7 and (1) of \S3.1. Using (1) of \S3.1, Lemma 2.2 one obtains
the following equations: 
$$
\align
& a_{2\Lambda_0}=a_{2\Lambda_1}=a_{2\Lambda_{l-1}}=a_{2\Lambda_l}=id,
\\
& a_{\Lambda_0+ \Lambda_1}= \sigma_{\tilde\Lambda_1} + 
\sigma_{\tilde\Lambda_{2l-1}},  
a_{\Lambda_i}= \sigma_{\tilde\Lambda_i} + \sigma_{\tilde\Lambda_{2l-i}},
2\leq i\leq l-2, \\
& a_{\Lambda_{l-1}+ \Lambda_l}= \sigma_{\tilde\Lambda_{l-1}} + 
\sigma_{\tilde\Lambda_{l+1}}. 
\endalign
$$ 
Let $W_0$ be the subspace with basis $\sigma_{\tilde\Lambda_{k}}, 0\leq
k\leq 2l-1$. Assume that 
$W=W_0\oplus W_1$. Then $W_1$ has dimension $2$. \par
Let us show that
for  
any irreducible twisted representation $x$ of ${\Cal A}_{Spin(M)_2}$
(these are representations which do not appear in (1) of \S3.1), 
$\widehat{j}x\neq x$.  Since by (1) of Lemma 2.7 and (1) of \S3.1 we have
$$ \dot S_{2\Lambda_0x} + \dot S_{2\Lambda_1x}=0,$$  so by (2) of \S2.2 
$$
Y_{2\Lambda_1x} = \omega_{2\Lambda_1x}/ \omega_{x} =-1.
$$
This implies that $2\Lambda_1x= \widehat{j}x\neq x$. Since $dim W_1=2$, 
and 
$$
\langle a_x, a_y\rangle = \langle x, y \rangle + \langle x,
\widehat{j} y \rangle
$$ by Lemma 2.2, it follows that 
there are exactly 4 twisted representations to generate 
two dimensional $W_1$.
But there are exactly 4 twisted representations in (4),  so these are all the
irreducible representations.\par
The univalences of the representations in (4) are given by 
(cf. \S1.4 of [KW]):
$$
\align
& \omega_{\widehat{1}}=\omega_{\widehat{j}}=1, 
\omega_{\widehat{\phi_k}} = \exp(\frac{\pi i k(l-k)}{4l}), 1\leq k
\leq l-1, \\
& \omega_{\widehat{\sigma_1}}= \omega_{\widehat{\sigma_2}}= 
\exp(\frac{\pi i (2l-1)}{8}),
\omega_{\widehat{\tau_1}}= 
\omega_{\widehat{\tau_2}}= -\exp(\frac{\pi i(2l-1) }{8}).
\endalign
$$
Using (2) of Lemma 2.7 the $\dot T$ matrix can be chosen to
be
$$
T_{xy} = \delta_{x.y} \omega_x \exp(\frac{-\pi i(2l-1)}{12}) \tag 6
$$    
where $\omega_x$ is given as above.
The $\dot S$-matrix is obtained by using  (1) to (7) of \S2.1 and Lemma
2.7,
and it is given in table 3.2. \par
\midinsert
\vskip4pt
\centerline{\vbox{\offinterlineskip\halign{
       \vrule # &\strut \ \hfil$#$\hfil \ &\vrule#
       &\strut \ \hfil$#$\hfil \
       &\vrule#
       &\strut \ \hfil$#$\hfil \
       &\vrule#
       &\strut \ \hfil$#$\hfil \
       &\vrule#
       &\strut \ \hfil$#$\hfil \
       &\vrule#
       &\strut \ \hfil$#$\hfil \
       &\vrule#
       &\strut \ \hfil$#$\hfil \
       &\vrule#\cr
      \noalign{\hrule}
      height1pt&\omit&&\omit&&\omit&&\omit&&\omit&&\omit&\cr  
      & \sqrt{8l}\times \dot S && \widehat{1}  && \widehat{j} && 
\widehat{\phi_l^j} && \widehat{\phi_{k'}} && \widehat{\sigma_j} && \widehat{\tau_j} & \cr
      height1pt&\omit&&\omit&&\omit&&\omit&&\omit&&\omit&\cr
        \noalign{\hrule}
      height1pt&\omit&&\omit&&\omit&&\omit&&\omit&&\omit&\cr
      & \widehat 1 && 1 && 1 && 1 && 2 && \sqrt{l} &&
      \sqrt{l} & \cr
      height1pt&\omit&&\omit&&\omit&&\omit&&\omit&&\omit&\cr
          \noalign{\hrule}
          \noalign{\hrule}
      height1pt&\omit&&\omit&&\omit&&\omit&&\omit&&\omit&\cr
      & \widehat{j} && 1 && 1 && 1 && 2 && -\sqrt{l} &&
      -\sqrt{l} & \cr  
          \noalign{\hrule}
      height1pt&\omit&&\omit&&\omit&&\omit&&\omit&&\omit&\cr
      & \widehat{\phi_l^i} && 1 && 1 && (-1)^l && 2(-1)^{k'} &&
       b_{ij}
       && b_{ij}  & \cr  
          \noalign{\hrule}
      height1pt&\omit&&\omit&&\omit&&\omit&&\omit&&\omit&\cr
      & \widehat{\phi_k} && 2 && 2 && 2(-1)^k &&
       4\cos{\frac{\pi kk'}{2l}} &&  0 && 0 & \cr  
          \noalign{\hrule}
      height1pt&\omit&&\omit&&\omit&&\omit&&\omit&&\omit&\cr
      & \widehat{\sigma_i} && \sqrt{l} && -\sqrt{l} &&
       b_{ij}  &&
       0 &&  a_{ij} &&
       -a_{ij} &\cr 
          \noalign{\hrule}
      height1pt&\omit&&\omit&&\omit&&\omit&&\omit&&\omit&\cr
      & \widehat{\tau_i} && \sqrt{l} && -\sqrt{l} && b_{ij}  &&
       0 &&  -a_{ij} && a_{ij} & \cr  
          \noalign{\hrule}
       }}}   
\vskip-10pt
\captionwidth{\hsize}
\botcaption{\bf Table 3.2}{\it $\sqrt{8l}\times \dot S$-matrix. Here
$a_{ij}= \sqrt{\frac{l}{2}}(1+(2\delta_{i,j}-1)\exp({-\frac{\pi il}{2}}))$, 
$b_{ij}=(-1)^{l+\delta_{i,j}}\sqrt{l}\exp{\frac{\pi il}{2}}$,
$\delta_{i,j}$
is the usual Delta function, and
$1\leq i,j\leq2$.}\endcaption
\vskip10pt
\endinsert             

Let us explain in more details how table 3.2 is obtained. First note that
if $\beta$ is a twisted representation, from (1) of Lemma 2.7 and
(1) we have
$$
\dot S_{\widehat{\phi_l^1} \beta} = - \dot S_{\widehat{\phi_l^2} \beta},
\dot S_{\widehat{1} \beta} = - \dot S_{\widehat{j} \beta},
\dot S_{\widehat{\phi_k} \beta}=0, 1\leq k\leq l-1. 
$$
Turning to nontwisted representations 
one obtains directly from (1) of Lemma 2.7 that
$$
\dot S_{\widehat{\phi_k} \widehat{\phi_k'}}= 4\cos(\frac{\pi k k'}{l}),
1\leq k, k'\leq l-1.
$$
By (1) of Lemma 2.7 we also have:
$$
\align
\dot S_{\widehat{\phi_l^1} \widehat{\phi_k'}} + 
\dot S_{\widehat{\phi_l^2} \widehat{\phi_k'}} 
&= \frac{(-1)^{k'}}{l}, 
\dot S_{\widehat{\phi_l^1} \widehat{\phi_l^2}} + 
\dot S_{\widehat{\phi_l^1} \widehat{\phi_l^1}} 
= \frac{(-1)^l}{2l}; \\
\dot S_{\widehat{\phi_l^1} \widehat{j}} + 
\dot S_{\widehat{\phi_l^2} \widehat{j}} 
&= \frac{1}{2l}, 
\dot S_{ \widehat{j}  \widehat{\phi_k'}} + 
\dot S_{ \widehat{1}  \widehat{\phi_k'}} = \frac{1}{2l}.
\endalign  
$$
These equations and (4) are sufficient to determine the rows 
labeled by $$\widehat{1}, 
\widehat{j}, \widehat{\phi_k}, \widehat{\phi_l^i}, 1\leq k\leq l-1,
i=1,2$$  
in Table 3.2. By the orthogonality of $\dot S$ matrix, 
the $\widehat{\phi_l^1}$-th row is orthogonal to the
$\widehat{1}$-th row, and we get the following 
equation
$ \dot S_{\widehat{\phi_l^i} \widehat{\sigma_1}} = - 
\dot S_{\widehat{\phi_l^i} \widehat{\sigma_2}}. 
$ 
This and (5) of \S2.1 are sufficient to determine all the entries in
table 3.2. \par
One checks easily that the genus $0$ modular matrices as given above
coincide with the the  modular matrices as given in Chap. 13 of [Kac].\par
Let us record the above results in the following:
\proclaim{Proposition 3.3}
All the irreducible covariant representations of ${\Cal A}_{Spin(M)_2}$
are listed in (4). The genus $0$ modular matrices of 
 ${\Cal A}_{Spin(M)_2}$ are given by (6) and table 3.2.
\endproclaim
\subheading{\S3.3. The diagonal coset 
$Spin(M)_2\subset Spin(M)_1\times Spin(M)_1$}  
In this section we determine the modular matrices for diagonal coset
$Spin(M)_2\subset Spin(M)_1\times Spin(M)_1$. The diagonal cosets
in the type $A$ case have been discussed in \S4.3 of [X2]. We will use the
same idea in \S4.3 of [X2] to obtain the  the modular matrices for diagonal 
coset $Spin(M)_2\subset Spin(M)_1\times Spin(M)_1$.\par
First we note that Th. 2.3 of [X2] holds in our case: the only change
one needs to make in the proof of Th. 2.3 of [X2] is to replace 
fermions by real ferminions as in [Bo].  
Since  $Spin(M)_2\subset Spin(2M)_1$ is cofinite by Lemma 3.2,
and  $Spin(M)_1\times Spin(M)_1\subset Spin(2M)_1$, by Prop. 3.1 of
[X2] $Spin(M)_2\subset Spin(M)_1\times Spin(M)_1$ is also cofinite. \par
Recall the set of integrable weights of 
irreducible positive energy representations of
$LSpin(2l)$ at level $k$ is given by:
$$
\align
P_+^{(k)}& = \{ \lambda= \lambda_0 \Lambda_0 + \lambda_1\Lambda_1 +...
\lambda_{l-1} \Lambda_{l-1} + \lambda_l \Lambda_l| \lambda_i\in {\Bbb N},\\
& \lambda_0 + \lambda_1 + 2( \lambda_2+...+\lambda_{l-2})+\lambda_{l-1}+
\lambda_l =k \}
\endalign
$$  
When $l$ is even, this set admits a ${\Bbb Z}_2 \times {\Bbb Z}_2$
automorphism
generated by 
$A_s, A_v$ where $A_s, A_v$ are given by:
$$
\align
& A_s (\lambda_0 \Lambda_0 + \lambda_1\Lambda_1 +...
\lambda_{l-1} \Lambda_{l-1} + \lambda_l \Lambda_l) = 
\lambda_0 \Lambda_l + \lambda_1\Lambda_{l-1} +...
\lambda_{l-1} \Lambda_{1} + \lambda_l \Lambda_0, \\ 
& A_v (\lambda_0 \Lambda_0 + \lambda_1\Lambda_1 +...
\lambda_{l-1} \Lambda_{l-1} + \lambda_l \Lambda_l) = 
\lambda_0 \Lambda_1 + \lambda_1\Lambda_0 +
\lambda_2\Lambda_2+ ... \lambda_{l-2} \Lambda_{l-2}
+ \\
& \lambda_{l-1} \Lambda_{l} + \lambda_l \Lambda_{l-1}. 
\endalign
$$ 
When $l$ is odd, this set admits a ${\Bbb Z}_4$
automorphism
generated by
$A_s$ where $A_s$ is given by
$$
\align
A_s (\lambda_0 \Lambda_0 + \lambda_1\Lambda_1 +...
\lambda_{l-1} \Lambda_{l-1} + \lambda_l \Lambda_l) & = 
\lambda_0 \Lambda_l + \lambda_1\Lambda_{l-1} +
\lambda_2\Lambda_{l-2}+ ... \lambda_{l-2} \Lambda_{2}
+ \\
& \lambda_{l-1} \Lambda_{0} + \lambda_l \Lambda_{1}. 
\endalign
$$
These automorphisms will be called {\it diagram automorphisms}. 
Note that the set of diagram automorphisms is 
${\Bbb Z}_2 \times {\Bbb Z}_2$ when $l$ is even, and 
${\Bbb Z}_4$ when $l$ is odd. \par
Th. 4.3 of [X2] now holds for our diagonal coset, with the action of
${\Bbb Z}_N$ there replaced by 
${\Bbb Z}_2 \times {\Bbb Z}_2$ when $l$ is even, and 
${\Bbb Z}_4$ when $l$ is odd, since the proof of [X2] applies verbatim.
As in \S4.3 of [X2], we denote by $[\dot \Lambda,\ddot \Lambda; \Lambda]$ 
the orbit of $(\dot \Lambda,\ddot \Lambda; \Lambda)$ under the diagonal
action of the diagram automorphisms. The the following is a 
list of irreducible representations of 
${\Cal A}_{Spin(M)_1\times Spin(M)_1/
Spin(M)_2}$.  
$$
\align
& 1:=[\dot \Lambda_0, \ddot \Lambda_0; \Lambda_0], 
j:= [\dot \Lambda_0, \ddot \Lambda_0; 2\Lambda_1], \\
& \phi_l^1:= [\dot \Lambda_0, \ddot \Lambda_0; 2\Lambda_{l-1}], 
\phi_l^2:= [\dot \Lambda_0, \ddot \Lambda_0; 2\Lambda_l], \text{\rm if}
\ \ l \in 2{\Bbb Z}, \\
& \phi_l^1:= [\dot \Lambda_0, \ddot \Lambda_1; 2\Lambda_{l-1}], 
\phi_l^2:= [\dot \Lambda_0, \ddot \Lambda_1; 2\Lambda_l], \text{\rm if}
\ \ l \in 2{\Bbb Z}+1, \\
& \phi_1 :=  [\dot \Lambda_0, \ddot \Lambda_1; \Lambda_0+\Lambda_1 ],  
\phi_2 :=  [\dot \Lambda_0, \ddot \Lambda_0; \Lambda_2 ], 
\phi_3 :=  [\dot \Lambda_0, \ddot \Lambda_1; \Lambda_3 ], ...,\\
& \phi_{l-2} :=  [\dot \Lambda_0, \ddot \Lambda_0; \Lambda_{l-2}]
, \text{\rm if}
\ \ l \in 2{\Bbb Z}, 
\phi_{l-2} :=  [\dot \Lambda_0, \ddot \Lambda_1; \Lambda_{l-2}]
, \text{\rm if}
\ \ l \in 2{\Bbb Z}+1, 
\\
& \phi_{l-1} :=  [\dot \Lambda_0, \ddot \Lambda_1; \Lambda_{l-1}+ \Lambda_{l}]
, \text{\rm if}
\ \ l \in 2{\Bbb Z}, 
\phi_{l-1} :=  [\dot \Lambda_0, \ddot \Lambda_0; \Lambda_{l-1}+ \Lambda_{l}]
, \text{\rm if}
\ \ l \in 2{\Bbb Z}+1, 
\\
& \sigma_1:=[\dot \Lambda_0, \ddot \Lambda_{l-1}; \Lambda_0+ \Lambda_{l-1}], 
\tau_1:=[\dot \Lambda_0, \ddot \Lambda_{l-1}; \Lambda_1+ \Lambda_{l}], 
\\
&\sigma_2:=[\dot \Lambda_0, \ddot \Lambda_{l}; \Lambda_0+ \Lambda_{l}], 
\tau_2:=[\dot \Lambda_0, \ddot \Lambda_{l}; \Lambda_1+ \Lambda_{l-1}] \tag 7 
\endalign
$$
The above is in fact a complete list of all the irreducible covariant
representations of ${\Cal A}_{Spin(M)_1\times Spin(M)_1/
Spin(M)_2}$ by Cor. 3.2 of [X4], since the proof of Cor. 3.2 of [X4]
applies verbatim in the present case.  We have also chosen our notations
to make the comparisons with the notations of [DVVV] easy (our $l$
corresponds to $N$ on P. 517, P. 518 of [DVVV]). \par
The univalences of the above representations are given by:
$$
\align
& \omega_1=\omega_j=1, \omega_{\phi_k} = \exp(\frac{\pi i k^2}{4l}), 1\leq k
\leq l-1, \\
& \omega_{\sigma_1}= \omega_{\sigma_2}= \exp(\frac{\pi i }{8}),
\omega_{\tau_1}= \omega_{\tau_2}= -\exp(\frac{\pi i }{8}) 
\endalign
$$
Using the remark after Prop. 3.1 of [X4] the $T$ matrix can be chosen to
be 
$$
T_{xy} = \delta_{x.y} \omega_x \exp(\frac{-\pi i}{12}) \tag 8
$$ where $\omega_x$ is given as above.
By Prop. 3.4 and (7), one can easily determine the $S$-matrix
(cf. (2) of Lemma 2.2 in [X3]) for 
${\Cal A}_{Spin(M)_1\times Spin(M)_1/Spin(M)_2}$. We have chosen our
notations in (7) so that the $S$-matrix is given by table 3.2 with
all the hats removed. We record the results of this section in the
following:
\proclaim{Proposition 3.4}
All the irreducible covariant representations of  
$$
{\Cal A}_{Spin(M)_1\times Spin(M)_1/Spin(M)_2}
$$ are given in (7), and its
modular matrices are given by (8) and table 3.2 ( with
all the hats removed).
\endproclaim
\subheading {\S3.4. The ${\Bbb Z}_2$ orbifold}
Note that the ${\Bbb Z}_2$ action on  ${\Cal A}_{U(1)_{2l}}$ given by $AdJ$
as defined before Lemma 3.1 is a proper action on ${\Cal A}_{U(1)_{2l}}$.
The reader familiar with  [FLM] may notice that the action given by $AdJ$
corresponds to $-1$ isometry of rank one lattice vertex operator algebras
(cf. [FLM] and [DN]). We have:
\proclaim{Lemma 3.5}
$$
{\Cal A}_{Spin(M)_1\times Spin(M)_1/Spin(M)_2}=
{\Cal A}_{U(1)_{2l}}^{{\Bbb Z}_2}.
$$
\endproclaim
\demo{Proof}
By definition we have a natural inclusion 
$$
{\Cal A}_{Spin(M)_1\times Spin(M)_1/Spin(M)_2}(I) \subset {\Cal A}_{U(1)_{2l}}^{{\Bbb Z}_2}(I), 
\forall I\in {\Cal I}
$$ since 
${\Cal A}_{Spin(M)_1\times Spin(M)_1}(I)$ is $AdJ$ invariant.
Note that
the inclusion has finite index $d$ since $Spin(M)_1\times Spin(M)_1\subset
Spin(2M)_1$
is a conformal inclusion which has finite index by [Bo]. 
By Prop. 2.4  we have
$$
\mu_{{\Cal A}_{{Spin(M)_1\times Spin(M)_1/Spin(M)_2}}} \times d^2 =
\mu_{{\Cal A}_{U(1)_{2l}}^{{\Bbb Z}_2}}. 
$$
But $\mu_{{\Cal A}_{Spin(M)_1\times Spin(M)_1/Spin(M)_2}}= 8l$ by Prop. 3.4, 
$\mu_{{\Cal A}_{U(1)_{2l}}^{{\Bbb Z}_2}}=8l$ by Th. 2.6, it follows that
$d^2=1$, and so ${\Cal A} = {\Cal A}_{U(1)_{2l}}^{{\Bbb Z}_2}$. \par
\enddemo
\qed \par
By Lemma 3.4, Th. 2.6 and  Prop. 3.4, we have proved the following theorem:
\proclaim{Theorem 3.6}
All the irreducible representations of ${\Cal A}_{U(1)_{2l}}^{{\Bbb Z}_2}$ 
are given by (7) for $l\geq 3$. 
These  irreducible representations give rise to a unitary
modular category whose genus $0$ 
modular matrices are given by (8) and table 3.2 with
all the hats removed.
\endproclaim
When $l=2$, $Spin(4)= SU(2)\times SU(2)$, one checks that Th. 3.6 
still holds in this case, where
the integrable weights of $LSpin(4)$ should be replaced by 
the integrable weights of $LSU(2)\times LSU(2)$.
 When $l=1$, Using the fact that
$ {\Cal A}_{U(1)_{2}}= {\Cal A}_{SU(2)_1}$ one can check easily
that 
$$
{\Cal A}_{U(1)_{2}}^{{\Bbb Z}_2}={\Cal A}_{U(1)_{8}}, 
$$
and ${\Cal A}_{U(1)_{8}}$ has already been studied in \S3.5 of [X5]. 
As noted before Th. 3.6, the 
${\Bbb Z}_2$ action on  ${\Cal A}_{U(1)_{2l}}$ given by $AdJ$
corresponding to $-1$ isometry of rank one lattice vertex operator algebras
(cf. [FLM] and [DN]). The classification of irreducible representations
of the orbifold rank one lattice vertex operator algebras is given in [DN]
which corresponds to the first part of Th. 3.6. We note that the $S$ matrix
can be identified with the $S$ matrix on P. 517 and P. 518 of [DVVV]. 
However there
are mistakes in the $S$ matrix on P. 517 and P. 518 of [DVVV]
corresponding to the entries of $a_{ij}, b_{ij}$ in Table 3.2. 
Table 3.2 gives the correct $S$ matrix. \par
\heading \S4. More examples and questions. \endheading
The lattice vertex operator algebras and their automorphim groups
provide a rich source of examples of orbifolds (cf. [FLM], [DN], [DGR]). 
We have determined the
genus $0$ modular matrices for the orbifold of rank $1$ lattice VOAs 
in Th. 3.6. 
It will be
interesting to generalize this to higher rank cases.\par 
Let ${\Cal A}_{SU(N)_k}$ be the irreducible conformal precosheaf and let $G$
be a finite subgroup of $SU(N)$. Then there is a natural action of
$G$ on ${\Cal A}_{SU(N)_k}$ and it is easy to check that this action is proper
(If the action of $G$ is not faithful, one can replace $G$ by a quotient
$G'$ as explained in the footnote of definition 2.1). Hence Th. 2.6 applies
in this case, and  we have a new family of unitary modular categories. 
It will be
interesting to study these modular categories. As a speical case, let
$N=2$. The finite subgroups of $SU(2)$ are classified into $A-D-E$ series.
Consider ${\Cal A}_{SU(2)_k}$ and $G/{\Bbb Z_2}= {\Bbb Z_k}$,
where ${\Bbb Z_2}$ is the center of $SU(2)$.  Note that 
the coset $$U(1)_{2k}\subset SU(2)_k$$ has been studied in \S3.5 of [X5].
We claim that 
$$ {\Cal A}_{SU(2)_k}^{G}= {\Cal A}_{U(1)_{2k}} \otimes
{\Cal A}_{SU(2)_k/U(1)_{2k}}.$$ 
Note that $G\subset U(1)$, and by definition we have
$$
 {\Cal A}_{U(1)_{2k}}(I) \otimes
{\Cal A}_{SU(2)_k/U(1)_{2k}}(I) \subset  {\Cal A}_{SU(2)_k}^{G}(I).
$$ 
But by Lemma 2.2 and \S3.5 of [X5], 
$$
\mu_{{\Cal A}_{SU(2)_{k}}} |k|^2= 
\mu_{{\Cal A}_{U(1)_{2k}}}\mu_{{\Cal A}_{SU(2)_k/U(1)_{2k}}}.
$$
It follows by Prop. 2.4 that the inclusion 
$$
 {\Cal A}_{U(1)_{2k}}(I) \otimes
{\Cal A}_{SU(2)_k/U(1)_{2k}}(I) \subset  {\Cal A}_{SU(2)_k}^{G}(I)
$$ has index $1$ which shows that
 $${\Cal A}_{SU(2)_k}^{G}= {\Cal A}_{U(1)_{2k}} \otimes
{\Cal A}_{SU(2)_k/U(1)_{2k}}.$$ 
Thus the unitary modular categories associated with ${\Cal
A}_{SU(2)_k}^{G}$ and the corresponding 3-manifold invariants
are determined by \S3.5 of [X5]. \par
Finally, let us mention that permutation orbifolds (cf. [BDM] and
references therein) provide another
interesting class of orbifolds. Let us formulate these orbifolds
in our setting. Let ${\Cal A}$ be an irreducible conformal 
precosheaf. Then ${\Cal A}$ tensor product itself $n$ times
${\Cal A}^{\otimes^n}:={\Cal A}\otimes {\Cal A}\otimes ... \otimes {\Cal A}$ 
is also an irreducible conformal 
precosheaf. Let $G\subset S_n$ be a finite subgroup of $S_n$, the 
permutation group on $n$ letters. Note that
any finite group is embedded in a permutation
group by Cayley's theorem. There is an obvious action of $G$ on
${\Cal A}^{\otimes^n}$ by permuting the $n$ tensors, and one checks
directly
by definitions that this action of $G$ on ${\Cal A}^{\otimes^n}$ is 
proper as defined in \S2.3. 
Note that if ${\Cal A}$ is $\mu$-rational, so is 
${\Cal A}^{\otimes^n}$ by definition. So if ${\Cal A}$ is $\mu$-rational,
by Th. 2.6 
we obtain a family of new unitary modular categories from the orbifold
$({\Cal A}^{\otimes^n})^G$. Consider a simple example where ${\Cal A}=
{\Cal A}_{SU(2)_1}$ and $n=2$. Let $G={\Bbb Z}_2=S_2$. Then by definition
one has ${\Cal A}_{SU(2)_1}\otimes {\Cal A}_{SU(2)_1}={\Cal
A}_{SU(2)_1\times
SU(2)_1}$, and  
$$
{\Cal A}_{SU(2)_2}(I) \otimes {{\Cal A}_{SU(2)_1\times
SU(2)_1/SU(2)_2}}(I) \subset ({\Cal A}^{\otimes^2})^G(I).
$$ 
By computing the $\mu$-index from \S3 of [X3]
we have
$$
\mu_{{\Cal A}_{SU(2)_2}} \mu_{{\Cal A}_{SU(2)_1\times
SU(2)_1/SU(2)_2}} = 4 \mu_{{\Cal A}_{SU(2)_1}}^2 
$$
and it follows from Prop. 2.4 as in the previous paragraph that
$$
{\Cal A}_{SU(2)_2} \otimes {{\Cal A}_{SU(2)_1\times
SU(2)_1/SU(2)_2}} = ({\Cal A}^{\otimes^2})^G.
$$ 
The modular category associated with $({\Cal A}^{\otimes^2})^G$  
and the corresponding 3-manifold invariants are therefore determined by
Prop. 3.6.3 of [X5].\par
In general, the modular matrices of the unitary modular 
categories associated with 
the permutation orbifolds above 
have been written down based on heuristic physics arguments in
[Ba]. It will be interesting to do the computations in our 
framework as in \S3 and compare with the results of [Ba].   
\heading References \endheading 
\roster 
\item"[Ba]" P. Bantay, {\it Permutation orbifolds}, hep-th/9910079.
\item"[BDM]" K. Barron, C. Dong and G. Mason, {\it Twisted sectors
for tensor product vertex operator algebras associated to permutation
groups.} math. QA/9803118
\item"[Bo]" J. B\"{o}ckenhauer, {\it An algebraic formulation of
level one Wess-Zumino-Witten models.} Rev. Math. Phys. {\bf 8} (1996) 925-948
\item"[BE1]" J. B\"{o}ckenhauer, D. E. Evans,
{\it Modular invariants, graphs and $\alpha$-induction for
nets of subfactors. I.},
Comm.Math.Phys., {\bf 197}, 361-386, 1998
\item"[BE2]" J. B\"{o}ckenhauer, D. E. Evans,
{\it Modular invariants, graphs and $\alpha$-induction for
nets of subfactors. II.},
Comm.Math.Phys., {\bf 200}, 57-103, 1999   
\item"[BE3]" J. B\"{o}ckenhauer, D. E. Evans,
{\it Modular invariants, graphs and $\alpha$-induction for
nets of subfactors. III.},
Comm.Math.Phys., {\bf 205}, 183-228, 1999
\item"[BE3]" J. B\"{o}ckenhauer, D. E. Evans,
{\it Modular invariants, graphs and $\alpha$-induction for
nets of subfactors. III.},
Comm.Math.Phys., {\bf 205}, 183-228, 1999. Also see hep-th/9812110 
\item"[BEK1]" J. B\"{o}ckenhauer, D. E. Evans, Y. Kawahigashi,
{\it On $\alpha$-induction, chiral generators and modular invariants
for subfactors}, Comm.Math.Phys., {\bf 208}, 429-487, 1999. Also
see math.OA/9904109
\item"[BEK2]" J. B\"{o}ckenhauer, D. E. Evans, Y. Kawahigashi,
{\it Chiral structure of modular invariants for subfactors},
Comm.Math.Phys., {\bf 210}, 733-784, 2000    
\item"{[B]}" R. Borcherds, {\it Vertex algebras, Kac-Moody
algebras and the Monster,}, Proc. Nat. Acad. Sci., U.S.A. 83, 3068-3071
(1986)
\item"[DM]" C. Dong and G. Mason, {\it On quantum Galois theory.} \par
Duke Math. J. {\bf 86}(1997), 305-321.
\item"[DVVV]" R. Dijkgraaf, C. Vafa, E. Verlinde and H. Verlinde,
{\it The operator algebra of orbifold models,} 
Comm.Math.Phys., {\bf 123}, 429-487, 1989.
\item"[DN]" C. Dong and K. Nagatomo, {\it Representations of vertex
operator
algebra $V_L^+$ for rank one lattice L.} math. QA/9807168
\item"[DGR]" C. Dong, R. Griess Jr. and A. Ryba, {\it 
Rank one lattice type vertex operator algebras and their automorphism
groups, II: E-series}. math.QA/9809022
\item"{[FJ]}" K. Fredenhagen and M. J\"{o}r$\beta$, {\it 
Conformal Haag-Kastler nets, pointlike localized fields and the
existence of opeator product expansions}. \par
Comm.Math.Phys., {\bf 176}, 541-554 (1996) 
\item"{[FLM]}" I. B. Frenkel, J. Lepowsky and J. Ries,
{\it Vertex operator algebras and the Monster}, Academic, New York, 1988. 
\item"{[FZ]}" I. Frenkel and Y. Zhu,
{\it Vertex operator algebras associated to representations of
affine and Virasoro algebras}, \par
Duke Math. Journal (1992), Vol. 66, No. 1, 123-168
\item"{[FG]}" J. Fr\"{o}hlich and F. Gabbiani, {\it Operator algebras and
Conformal field theory, } Comm. Math. Phys., {\bf 155}, 569-640 (1993). 
\item"{[GL]}"  D.Guido and R.Longo, {\it  The Conformal Spin and
Statistics Theorem},  \par
Comm.Math.Phys., {\bf 181}, 11-35 (1996)  
\item"{[ILP]}" M. Izumi, R. Longo and S. Popa, {\it 
A Galois correspondence for compact groups of automorphisms of
von Neumann Algebras with a generalization to Kac algebras.}
J. Funct. Analysis, {\bf 155}, 25-63 (1998) 
\item"{[Kac]}"  V. G. Kac, {\it Infinite Dimensional Lie Algebras}, 3rd
Edition,
Cambridge University Press, 1990.
\item"{[KLM]}" Y. Kawahigashi, R. Longo and M. M\"{u}ger,
{\it Multi-interval Subfactors and Modularity of Representations in
Conformal Field theory}, Preprint 1999, see also math.OA/9903104.
\item"{[KW]}"  V. G. Kac and M. Wakimoto, {\it Modular and conformal
invariance constraints in representation theory of affine algebras},
Advances in Math., {\bf 70}, 156-234 (1988).	
\item"{[L1]}"  R. Longo, {\it Minimal index and braided subfactors}, J.
Funct. Anal., {\bf 109}, 98-112 (1992).
\item"{[L2]}"  R. Longo, {\it Duality for Hopf algebras and for subfactors},
I, Comm. Math. Phys., {\bf 159}, 133-150 (1994).
\item"{[L3]}"  R. Longo, {\it Index of subfactors and statistics of
quantum fields}, I, Comm. Math. Phys., {\bf 126}, 217-247 (1989.
\item"{[L4]}"  R. Longo, {\it Index of subfactors and statistics of
quantum fields}, II, Comm. Math. Phys., {\bf 130}, 285-309 (1990).
\item"{[L5]}"  R. Longo, Proceedings of International Congress of
Mathematicians, 1281-1291 (1994).   
\item"{[LR]}"  R. Longo and K.-H. Rehren, {\it Nets of subfactors},
Rev. Math. Phys., {\bf 7}, 567-597 (1995).  
\item"{[MS]}" G. Moore and N. Seiberg, {\it Taming the conformal zoo},
Lett. Phys. B  {\bf 220}, 422-430, (1989).  
\item"[PS]" A. Pressley and G. Segal, {\it Loop Groups,} O.U.P. 1986.
\item"[Reh]" Karl-Henning Rehren, {\it Braid group statistics and their
superselection rules} In : The algebraic theory of superselection
sectors. World Scientific 1990
\item"[T]"  M. Takesaki, {\it Conditional expectation in von Neumann
algebra,} J. Funct. Analysis 9 (1972), 306-321.
\item"{[TL]}" V. Toledano Laredo, {\it Fusion of
Positive Energy Representations of $LSpin_{2n}$}.
Ph.D. dissertation, University of Cambridge, 1997
\item"[Tu]" V. G. Turaev, {\it Quantum invariants of knots and
3-manifolds,} Walter de Gruyter, Berlin, New York 1994
\item"{[W1]}"  A. Wassermann, Proceedings of International Congress of
Mathematicians, 966-979 (1994)
\item"{[W2]}"  A. Wassermann, {\it Operator algebras and Conformal
field theories III},  Invent. Math. Vol. 133, 467-539 (1998)
\item"{[W3]}"  A. Wassermann, {\it  Operator algebras and Conformal
field theories}, preliminary notes of 1992
\item"{[W4]}"  A. Wassermann, with contributions by V. Jones,{\it
Lectures on operator algebras and conformal field theory},
Proceedings of Borel Seminar, Bern 1994, to appear 
\item"[X1]" F.Xu, {\it   New braided endomorphisms from conformal
inclusions, } \par
Comm.Math.Phys. 192 (1998) 349-403.
\item"[X2]" F.Xu, {\it Algebraic coset conformal field theories},
 math.OA/9810035, to appear in Comm. Math. Phys.
\item"[X3]" F.Xu, {\it Algebraic coset conformal field theories II}, \par
 math.OA/9903096, Publ. RIMS, vol.35 (1999), 795-824.
\item"[X4]" F.Xu, {\it On a conjecture of Kac-Wakimoto}, \par
 math.RT/9904098,
\item"[X5]" F.Xu, {\it 3-manifold invariants from cosets},
math.GT/9907077.    
\item"[X6]" F.Xu, {\it Jones-Wassermann subfactors for
Disconnected Intervals}, \par
 q-alg/9704003, to appear in Comm. Contemp. Math.
\endroster

\enddocument